\documentclass[12
pt,twoside]{amsart}
\usepackage{amssymb}
\usepackage{a4wide}
\usepackage[latin1]{inputenc}
\usepackage[T1]{fontenc}
\usepackage{times}

\def\hpq0{h^{p,q}_{\leq 0}}
\def\Hpq0{\H_{\leq 0}^{p,q}}
\def\gr{\otimes}

\def\dbar{\bar\partial}
\def\ddbar{\partial\dbar}

\def\R{{\mathbb R}}

\def\C{{\mathbb C}}

\def\H{{\mathcal H}}

\def\Re{{\rm Re\,  }}

\def\be{\begin{equation}}
\def\ee{\end{equation}}

\newtheorem{thm}{Theorem}[section]
\newtheorem{lma}[thm]{Lemma}

\newtheorem{prop}[thm]{Proposition}

\theoremstyle{definition}

\theoremstyle{remark}

\newtheorem{preremark}{Remark}
\newtheorem{preex}{Example}

\numberwithin{equation}{section}

\begin{document}

\title[]
{Curvature of vector bundles associated to holomorphic fibrations.}

\author[]{ Bo Berndtsson}

\address{B Berndtsson :Department of Mathematics\\Chalmers University
  of Technology 
  and the University of G\"oteborg\\S-412 96 G\"OTEBORG\\SWEDEN,\\}

\email{ bob@math.chalmers.se}

\begin{abstract}
{Let $L$ be a (semi)-positive line bundle over a Kähler manifold, $X$,
  fibered over a complex manifold $Y$. Assuming the fibers are compact
  and non-singular we prove that the hermitian vector bundle $E$ over
  $Y$ whose fibers over points $y$ 
  are the spaces of global sections over $X_y$ to $L\gr K_{X/Y}$
  endowed with the 
  $L^2$-metric is (semi)-positive in the sense of Nakano. We also
  discuss various
  applications, among them  a partial result on a conjecture of Griffiths
  on the positivity of ample bundles.  This is a revised and much
  expanded version of a previous preprint 
  with the title `` Bergman kernels and the curvature of vector bundles''.}
\end{abstract}

\bigskip

\maketitle

\section{Introduction}

Let us first consider a domain $D=U\times\Omega$  in $\C^m\times\C^n$ and a 
function $\phi$, plurisubharmonic in $D$. We also assume
for simplicity that $\phi$ is smooth up to the boundary . Then, for
each $t$ in $U$, 
$\phi^t(\cdot):=\phi(t,\cdot)$ is plurisubharmonic in $\Omega$ and we denote
by $A^2_t$ the Bergman spaces of holomorphic functions in $\Omega$
with norm
$$
\|h\|^2=\|h\|^2_t=\int_\Omega |h|^2e^{-\phi^t}. 
$$
The spaces $A^2_t$ are then all equal as vector spaces but have norms
that vary with $t$. The - infinite rank - vector bundle $E$ over $U$
with fiber $E_t=A^2_t$ is therefore trivial as a bundle but is
equipped with a nontrivial metric. The first result of this paper is
the following theorem. 
\begin{thm} If $\phi$ is (strictly) plurisubharmonic, then the
  hermitian bundle $(E, \|\cdot\|_t)$ is (strictly) positive 
  in the sense of Nakano.
\end{thm}

Of the two main differential geometric notions of positivity (see
section 2, where these matters are reviewed in the slightly non
standard setting of bundles of infinite rank), positivity in the sense
of Nakano is the stronger one and 
 implies the  weaker property of positivity in the sense of
Griffiths. On the other hand the Griffiths notion of positivity has
nicer functorial properties and implies in particular that the dual
bundle is negative (in the sense of Griffiths). This latter property
is in turn equivalent to the condition that if $\xi$ is  any
nonvanishing local
holomorphic section of the dual bundle, then  the function
$$
\log\|\xi\|^2
$$
is  plurisubharmonic. We can obtain such
holomorphic sections to the dual bundle from point evaluations. More 
precisely, let $f$ be a holomorphic map from $U$ to $\Omega$ and
define $\xi_t$ by is action on a local section of $E$  
$$ 
\langle \xi_t, h_t\rangle= h_t(f(t)).
$$
Since the right hand side here is a holomorphic function of $t$, $\xi$
is indeed a holomorphic section of $E^*$. The norm of $\xi$ at a point
is given by
$$
\|\xi_t\|^2=\sup_{\|h_t\|\leq 1}| h_t(f(t))|^2=K_t(f(t),f(t)),
$$
where $K_t(z,z)$ is the Bergman kernel function for $A^2_t$. It
therefore follows from Theorem 1.1 that $\log K_t(z,z)$ is
plurisubharmonic (or identically equal to $-\infty$) 
in $D$, which is the starting point of the  results in
\cite{Berndtsson}. The point here is of course that   $\log K_t(z,z)$
is plurisubharmonic with respect to the parameter $t$, and even with
respect to all the variables $(t,z)$ jointly. 

In \cite{Berndtsson} it is proved that this subharmonicity property of
the Bergman kernel persists if $D$ is a general pseudoconvex domain in
$\C^m\times\C^n$, for general plurisubharmonic weight functions. In
this case the spaces $A^2_t$ are the Bergman 
spaces for the slices of $D$, $D_t=\{z; (t,z) \quad \text{lies in}
\quad D\}$. This
more general case should also lead to a positively curved vector
bundle. The main problem in proving such an extension of Theorem 1.1 is not to prove the
inequalities involved, but rather to define the right notion of vector
bundle in this case. In general, the spaces $A^2_t$ will not be
identical as vector spaces, so the bundle in question is not locally
trivial.

There is however  a natural analog of Theorem 1.1 for holomorphic 
fibrations with compact fibers. Consider a complex manifold $X$ of
dimension $n+m$ which is smoothly fibered
over 
another connected complex $m$-dimensional manifold $Y$. There is then 
a holomorphic map, $p$, from 
$X$ to $Y$ with surjective differential,  and all the fibers
$X_t=p^{-1}(t)$ are assumed compact.  This implies, see
\cite{Voisin}, that the fibers are all diffeomorphic, but they are in
general not  biholomorphic to each other. We next need a substitute
for the assumption on pseudoconvexity in Theorem 1.1. What first comes
to mind is that $X$ be a projective fibration, i e that there be a
strictly positive line bundle on $X$. This would mean that $X$
contains a divisor $A$ such that $X\setminus A$ is Stein, so that we
would be almost in the Stein case, which is quite similar to the case in
Theorem 1.1. It turns out however that all we need to assume is that
$X$ be K\"ahler.

 Let  $L$ be a
holomorphic, hermitian line bundle over the total space $X$. Our
substitute for the Bergman spaces $A^2_t$ is now the space of global
sections over each fiber to $L\gr K_{X_t}$, 
$$
E_t=\Gamma(X_t,L|X_t\gr K_{X_t}),
$$
where $K_{X_t}$ is the canonical bundle of, i e the bundle of forms of
bidegree $(n,0)$ on, each fiber. We assume that $L$ is semipositive so
that the hermitian metric on $L$  has
nonnegative curvature form. Fix a point $y$ in $Y$ and choose local
coordinates $t=(t_1, ...t_m)$ near $y$ with $t(y)=0$. We consider the
coordinates as functions on $X$ by identifying $t$ with $t\circ p$,
and let $dt=dt_1\wedge ...dt_m$. The canonical bundle of a fiber $X_t$
can then be identified with the restriction of $K_X$, the canonical
bundle of the total space, to $X_t$ by mapping a local section $u$ to
$K_{X_t}$ to  $u\wedge dt$. This map is clearly injective, and it is
also surjective since any section $w$ to $K_X$ can locally be written
$$
w=u\wedge dt,
$$
and the restriction of $u$ to $X_t$ is independent of the choice of
$u$. With this identification, any global holomorphic section of $L\gr K_{X_t}$
over a fiber can be holomorphically extended to a holomorphic section of  
$L\gr K_{X_s}$, for $s$ near $t$. When $L$ is trivial this follows from the
Kähler assumption, by invariance of Hodge numbers, see
\cite{Voisin}. When $L$ is semipositive it follows from a variant of
the Ohsawa-Takegoshi extension theorem, that we discuss in an appendix
. Starting from a basis for 
$\Gamma(X_t,L|X_t\gr K_{X_t})$ we therefore get a local holomorphic
frame for $E$, so $E$ has a natural structure as a holomorphic vector
bundle. Moreover,
elements of
$E_t$ can be naturally integrated over the fiber and we obtain in this
way a metric, $\|\cdot\|$ on $E$ in complete analogy with the plane
case. We then get the same conclusion as before:
\begin{thm}If the total space $X$ is Kähler and  $L$ is (semi)positive
  over $X$, then  $(E,\|\cdot\|)$ is 
  (semi)positive in the sense of Nakano. 
\end{thm}
  
This can be compared to  results of Fujita, \cite{Fujita}, Kawamata,
\cite{Kawamata} and  Koll\'ar,
\cite{Koll\'ar},  who proved  positivity properties for $E$
when $L$ is 
trivial and $X$ projective. Kawamata also extended these results to
multiples of the canonical bundle. Related work is also due to Tsuji,
see \cite{Tsuji} and the 
references therein. The proofs in these papers are based on results of
Griffiths on variations of Hodge structures, whereas our proofs use
techniques related to Hörmander-like $L^2$-estimates for $\dbar$. This
seems to be the main novelty in our approach. Among the advantages are
that this permits us to work directly in the twisted context (i e with
nontrivial $L$), it gives Nakano positivity and not just Griffiths
positivity and it also works in the noncompact case (like Theorem
1.1). Moreover it gives explicit interesting lower bound for the
curvature operator; see section 6.

 After a first version of this paper was posted on the
ArXiv, Tsuji also announced a version of Theorem 1.2 in
\cite{Tsuji2}. In this paper he indicates a proof  that $E$ is
Griffiths positive, assuming the fibrations is projective, by a
 reduction to the case of a locally trivial fibration.

Theorems 1.1 and 1.2 have applications, or at least illustrations, in
a number of different contexts. One concerns subharmonicity
properties of  Bergman kernels depending on a parameter, which follow
from Theorem 1.1 as 
explained above, see 
\cite{Berndtsson}. This can be seen as a complex version of the 
Brunn-Minkowski inequality. The link to Brunn-Minkowski theory lies in the
fact that the inverse of the volume of a domain in $\R^n$ is the
Bergman kernel for the space of constant functions on the domain. The
classical Brunn-Minkowski inequality is therefore a convexity
property, with respect to parameters, 
of the Bergman kernel for the nullspace of the $d$-operator, whereas
here we have plurisubharmonicity of   the Bergman kernel for the
nullspace of the $\dbar$-operator. 

Just like in the case of Theorem 1.1, Theorem 1.2 also has as a
consequence a result on plurisubharmonicity of a Bergman kernel. In
this case however, the Bergman kernel is not a function, but
transforms as a metric on the {\it relative canonical bundle} of the
fibration
$$
K_{X/Y}:= K_X\gr  p^*(K_Y)^{-1},
$$
twisted with $L$. 
In  forthcoming work with Mihai Paun, \cite{B-Paun1} and
\cite{B-Paun2}, we show how Theorem 1.2 
implies that this Bergman kernel  
has semipositive curvature or is identically equal to 0. In particular
$L+ K_{X/Y}$ is 
pseudoeffective if it has a nontrivial square integrable holomorphic
section over 
at least one fiber. We also
extend this result to the case when $L$ has a 
singular metric and the map $p$ is not necessarily a smooth
fibration. Moreover,  these methods can be extended to twisted
multiples of the canonical bundle, generalising to the twisted case
the work of Kawamata quoted above.

Not surprisingly, the curvature of the bundle $E$ in Theorem 1.2 can
be zero at some point and in some direction only if the curvature of
the line bundle $L$ also degenerates. In section 5 we  prove a
result that indicates that conditions for degeneracy of the curvature
of $E$ are much more restrictive than that: When $X$ is a product,
null vectors for the curvature can only come from infinitesimal
automorphisms of the fiber. 

In section 6 we discuss some, largely philosophical, relations between
Theorem 1.2 and recent work on the variation of Kähler metrics. This
corresponds to the case when $X$ is a product $U\times Z$ with
one-dimensional base $U$, and when $L$ is the pull-back of a bundle on
$Z$ under the second projection map. The variation of the metric on
$L$ that we get from the fibration then gives a path in the space of
Kähler metrics on $Z$ and the lower bound that we get for the
curvature operator in this case is precisely the Toepliz operator
defined by the geodesic curvature of this path. This theme is further
developed in \cite{3berndtsson}.
 
Another example of the situation in Theorem 1.2   arises naturally if we
start with a (finite rank) holomorphic vector bundle $V$ over $Y$ and
let $\mathbb P(V)$ be the associated bundle of projective spaces  of
the dual bundle $V^*$. This is then clearly an - even  locally
trivial -
holomorphic fibration  and there is a naturally defined  line bundle
$L$ over the total space
$$
L=O_{\mathbb P(V)}(1),
$$
that restricts to the hyperplane section bundle over each fiber. The
global holomorphic 
sections of 
this bundle over each fiber are now the linear forms on $V^*$, i e the
elements of $V$. To
obtain sections to our bundle $E$ defined above, we  take tensor products
with the canonical bundle. We therefore replace $L$ by 
$$
L^{r+1}=O_{\mathbb P(V)}(r+1)
$$
(with $r$ being the rank of $V$). Since the canonical bundle of a
fiber is $O(-r)$ we see that on each fiber $L\gr K_{X_t}=O(1)$ so its
space of global sections is again equal to $V_t$. Define as before
$$
E_t=\Gamma(X_t,L^{r+1}|X_t\gr K_{X_t}).
$$
One can then verify that, globally,  $E$ is isomorphic to $V\gr\det V$. The
condition that $L$ is positive is now equivalent to $O_{\mathbb P(V)}(1)$ being
positive which is the same as saying that $V$ is {\it ample} in the
sense of Hartshorne, \cite{Harts}. We therefore obtain (in section 7)
the following result 
 as a corollary of Theorem 1.2.

\begin{thm} Let $V$ be a (finite rank) holomorphic vector bundle over
  a complex manifold which is ample in the sense of Hartshorne. Then
  $V\gr\det V$ has a smooth hermitian metric which is strictly
  positive in the sense of Nakano. 
\end{thm}
Replacing $O_{\mathbb P(V)}(r+1)$ by $O_{\mathbb P(V)}(r+m)$, we also
get that $S^m(V)\gr \det V$ is Nakano-positive for any non-negative
$m$, where $S^m(V)$ is the $m$:th symmetric power of $V$.

It is a well known conjecture of Griffiths, \cite{Griffiths}, that an ample
vector bundle is positive in the sense of Griffiths.
Theorem 1.3 can perhaps be seen as indirect evidence for this conjecture, since by a a theorem of
Demailly, \cite{Demailly 2}, $V\gr\det V$ is Nakano positive if $V$ itself is
Griffiths positive. It seems  that not so much is known about
Griffiths' conjecture in general, except that it does hold when $Y$ is a
compact curve (see \cite{Umeu}, \cite{Camp}).

After the first version of this mansucript was completed I received a preprint
 by C Mourougane and S Takayama, \cite{M}. There they  prove that
 $V\gr\det V$ is positive in 
the sense of Griffiths, assuming the base manifold is projective. The
method of proof  is  quite different from
this paper, as is the metric they find.

We end this introduction with a brief discussion of the proofs.
The proof of Theorem 1.1 is based on regarding the bundle $E$ as a
holomorphic subbundle of the hermitian bundle $F$ with fibers
$$
F_t=L^2(\Omega,e^{-\phi^t})=:L^2_t.
$$
 By
definition, the curvature of $F$ is a $(1,1)$-form 
$$
\sum \Theta^F_{j k} dt_j\wedge d\bar t_k
$$
whose coefficients are operators on $F_t$. By direct and simple
computation,
$$
\Theta^F
$$
is the operator of multiplication with
$\partial_t\bar\partial_t\phi$, so it is positive as soon as $\phi$
is plurisubharmonic of $t$ for $z$ fixed.  By a formula of Griffiths, the
curvature of the holomorphic subbundle $E$ is obtained from the
curvature of $F$ by subtracting the {\it second fundamental form } of
$E$, and the crux of the proof is to control this term by the
curvature of $F$. For this we note that the second fundamental form is
given by the square of the norm of an element in the orthogonal
complement of $A^2_t$ in $L^2_t$. This element is therefore the
minimal solution of a certain $\dbar$-equation, and the needed
inequality follows from an application of Hörmander's $L^2$-estimate.

We have not been able to generalize this proof to the situation of
Theorem 1.2. The proof does generalize to the case of a holomorphically
trivial fibration, but in the general case we have not been able to
find a natural complex structure on the space of all (not necessarily
holomorphic) $(n,0)$ forms, extending the complex structure on $E$. We
therefore compute directly the 
Chern connection of the bundle $E$ itself, and compute the curvature
from there, much as one  proves Griffiths' formula. In these
computations appears also the Kodaira-Spencer class of the fibration,
\cite{Voisin}. This class plays somewhat the role of another second
fundamental form, but this time of a quotient bundle, arising when we
restrict $(n,0)$-forms to the fiber. The Kodaira-Spencer class therefore
turns out to give a positive contribution to the curvature. This proof
could also be adapted to give Theorem 1.1 by using fiberwise complete
Kähler metrics, but we have chosen not to do so since the first proof
seems conceptually clearer.

Finally, I would like to thank Sebastien Boucksom for pointing out the
relation between Theorem 1.1 and the Griffiths conjecture, Jean-Pierre
Demailly for encouraging me to treat also  the case of a general
non-trivial fibration
and Yum-Tong Siu and Mihai Paun 
for helpful discussions. Last but not
least, thanks are due to Hiroshi Yamaguchi, whose work on 
plurisubharmonicity of the Robin function \cite{Yamaguchi} and Bergman
kernels, \cite{Maitani-Yamaguchi} was an important source of
inspiration for this work.

\section{Curvature of finite and infinite rank bundles}

Let $E$ be a holomorphic vector bundle with a hermitian metric over a
complex manifold $Y$. By 
definition this means that there is a holomorphic projection map $p$
from $E$ to $Y$ and that every point in $Y$ has a neighbourhood $U$
such that $p^{-1}(U)$ is isomorphic to $U\times W$, where $W$ is a
vector space equipped with a smoothly varying hermitian metric. In our
applications it is important to be able to allow this vector space to
have infinite dimension, in which case we assume that the metrics are
also complete, so that the fibers are Hilbert spaces. 

Let $t=(t_1,...t_m)$ be a system of local coordinates on $Y$. The
Chern connection, $D_{t_j}$ is now given by a collection of
differential operators acting on smooth sections to $U\times W$ and
satisfying
$$
\partial_{t_j}( u,v)=(D_{t_j}u,v) +( u,\dbar_{t_j}v),
$$ 
with $\partial_{t_j}=\partial/\partial_{t_j}$ and
$\dbar_{t_j}=\partial/\partial\bar t_j$. The curvature of the Chern 
connection is a $(1,1)$-form of operators 
$$
\Theta^E=\sum\Theta^E_{j k}dt_j\wedge d\bar t_k,
$$
where the coefficients $\Theta^E_{j k}$ are densily defined operators on
$W$. By definition these coefficients are the commutators
$$
\Theta^E_{j k}= [ D_{t_j},\dbar_{t_k}].
$$
The vector bundle is said to be positive in the sense of Griffiths if
for any section $u$  to $W$ and any vector $v$ in $\C^m$
$$
\sum(\Theta^E_{j k}u,u)v_j\bar v_k\geq \delta\|u\|^2|v|^2
$$
for some positive $\delta$. $E$ is said to be positive in the sense of
Nakano if for any $m$-tuple $(u_1,...u_m)$ of sections to $W$
$$
\sum(\Theta^E_{j k}u_j,u_k)\geq \delta\sum\|u_j\|^2
$$
Taking $u_j=u v_j$ we see that Nakano positivity implies positivity in
the sense of Griffiths. 

The dual bundle of $E$ is the vector bundle $E^*$ whose fiber at a
point $t$ in $Y$ is the Hilbert space dual of $E_t$. There is
therefore a natural  antilinear isometry
between $E^*$ and $E$, which we will denote by  $J$. If $u$ is a local
section of $E$, $\xi$ is 
a local section of $E^*$, and $\langle\cdot,\cdot\rangle$ denotes the
pairing between $E^*$ and $E$ we have
$$
\langle\xi,u\rangle=(u,J\xi).
$$
Under the natural holomorphic structure on $E^*$ we then have
$$
\dbar_{t_j}\xi=J^{-1}D_{t_j}J\xi,
$$
and the Chern connection on $E^*$ is given by
$$
D^*_{t_j}\xi=J^{-1}\dbar_{t_j}J\xi.
$$
It follows that
$$
\dbar_{t_j}\langle\xi,u\rangle=
\langle\dbar_{t_j}\xi,u\rangle +\langle\xi,
\dbar_{t_j}u\rangle,
$$
and 
$$
\partial_{t_j}\langle\xi,u\rangle=
\langle D^*_{t_j}\xi,u\rangle +\langle\xi,D_{t_j}u\rangle,
$$
and hence 
$$
0=\left[\partial_{t_j},\dbar_{t_k}\right]\langle\xi,u\rangle =
\langle\Theta^{E^*}_{j k}\xi,u\rangle +\langle\xi,\Theta^E_{j k}u\rangle,
$$ 
if we let $\Theta^{E^*}$ be the curvature of $E^*$.
If $\xi_j$ is an $r$-tuple of sections to $E^*$, and $u_j=J\xi_j$,  we
thus see that 
$$
\sum(\Theta^{E^*}_{j k}\xi_j,\xi_k)=-\sum(\Theta^E_{j k}u_k,u_j).
$$
Notice that the order between $u_k$ and $u_j$ in the right hand side is
opposite to the order between the $\xi$s in the left hand
side. Therefore $E^*$ is negative in the sense of Griffiths iff $E$ is
positive in the sense of Griffiths, but we can not draw the same
conclusion in the the case of Nakano positivity.

If $u$ is a holomorphic section of $E$ we also find that
$$
\frac{\partial^2}{\partial t_j\partial\bar t_k}(u,u)=
(D_{t_j}u,D_{t_k}u) -(\Theta^E_{j k}u,u)
$$
and it follows after a short computation that $E$ is (strictly)
negative in the 
sense of Griffiths if and only if $\log\|u\|^2$ is (strictly)
plurisubharmonic for any nonvanishing holomorphic section $u$.

We next briefly recapitulate the Griffiths formula for the curvature
of a subbundle. Assume $E$ is a holomorphic subbundle of the bundle
$F$, and let $\pi$ be the fiberwise orthogonal projection from $F$ to
$E$. We also let $\pi_\bot$ be the orthogonal projection on the
orthogonal complement of $E$.
By the definition of Chern connection we have 
$$
D^E=\pi D^F.
$$
Let $\dbar_{t_j}\pi$ be defined by
\begin{equation}
\dbar_{t_j}(\pi u)=(\dbar_{t_j}\pi)u +\pi(\dbar_{t_j} u).
\end{equation}
Computing the commutators occuring in the  definition of curvature we see that
\begin{equation}
\Theta^E_{j k} u=-(\dbar_{t_k}\pi)D^F_{t_j} u+\pi\Theta^F_{j k} u,
\end{equation}
if $u$ is a section of $E$. 
By (2.1) $(\dbar\pi)v=0$ if $v$ is a  section of $E$,
so
\begin{equation}
(\dbar\pi)D^F u=(\dbar\pi)\pi_\bot D^F u.
\end{equation}
 Since
$\pi \pi_\bot=0$ it also follows that
$$
(\dbar\pi)\pi_\bot D^F u=-\pi\dbar(\pi_\bot D^Fu).
$$
Hence,  if $v$ is also a section of $E$,
$$
((\dbar_{t_k}\pi)D^F_{t_j} u,v)=-(\dbar_{t_k}(\pi_\bot D^F_{t_j}u),v)=
$$
$$
=
((\pi_\bot
D^F_{t_j}u),D^F_{t_k}v)=(\pi_\bot(D^F_{t_j}u),\pi_\bot(D^F_{t_k}v)) .
$$
Combining with (2.2) we finally get that if $u$ and $v$ are both
sections to $E$ 
then 
\begin{equation}
(\Theta^F_{j k}u,v)=(\pi_\bot(D^F_{t_j}u),\pi_\bot(D^F_{t_k}v))+
(\Theta^E_{j k}u,v),
\end{equation}
and thus
$$
\sum (\Theta^F_{j k}u_j,u_k)= \|\pi_{\bot}\sum D_{t_j}^F u_j\|^2
+\sum(\Theta^E_{j k}u_j,u_k).
$$ 
which is the starting point for the proof of Theorem 1.1.

For the proof of Theorem 1.2 we finally describe another way of
computing the curvature form of a vector bundle. Fix a point $y$ in
$Y$ and choose local coordinates $t$ centered at $y$. Any point $u_0$ in
the fiber $E_0$ over $y$  can be extended to a  holomorphic section
$u$ of
$E$ near 0. Modifying $u$ by a linear combination $\sum t_j v_j$ for
suitably chosen local holomorphic sections $v_j$ we can also arrange
things so that $Du=0$ at $t=0$. Let $u$ and $v$ be two local sections
with this property and compute
$$
\partial_{\bar t_k}\partial_{ t_j}(u,v)= \partial_{\bar t_k}(D_{t_j}
u,v)=(\partial_{\bar t_k}D_{t_j}u,v)= -(\Theta^E_{j k}u,v).
$$
Let $u_j$ me an $m$-tuple of holomorphic sections to $E$, satisfying
$D u_j=0$ at 0. Put

$$
T_{u}=\sum (u_j,u_k)\widehat{dt_j\wedge d\bar t_k}.
$$
Here $\widehat{dt_j\wedge d\bar t_k}$ denotes the wedge product of all
$dt_i$ and $d\bar t_i$ except $dt_j$ and $d\bar t_k$, multiplied by a
constant of absolute value 1,
chosen so that $T_u$ is a positive form.
Then
\be
i\ddbar T_{u}=-\sum(\Theta^E_{j k}u_j,u_k)dV_t,
\ee
so  $E$ is Nakano-positive at a given point if and only if this expression
is negative for any choice of holomorphic sections $u_j$ satisfying
$Du_j=0$ at the point.

\section{ The proof of Theorem 1.1}

We consider the setup described before the statement of Theorem 1.1
in the introduction. Thus $E$ is the vector bundle over $U$ whose
fibers are the Bergman spaces $A^2_t$ equipped with the weighted $L^2$
metrics induced by  $L^2(\Omega,e^{-\phi^t})$. We also let $F$ be the
vector bundle with fiber  $L^2(\Omega,e^{-\phi^t})$, so that $E$ is a
trivial subbundle of the trivial bundle $F$ with a  metric
induced from a nontrivial metric on $F$. From the definition of the
Chern connection we see that
$$
D^F_{t_j}=\partial_{t_j}-\phi_j,
$$
where the last term in the right hand side should be interpreted as the operator of
multiplication by the (smooth) function
$-\phi_j=-\partial_{t_j}\phi^t$. (In the sequel we use the letters
$j, k$ for indices of the $t$-variables,and the letters
$\lambda,\mu$ for indices of the $z$-variables.) For the
curvature of $F$ we therefore get
$$
\Theta^F_{j k}=\phi_{j k},
$$
the operator of multiplication with the complex Hessian of $\phi$ with respect
to the $t$-variables. We shall now apply formula (2.4), so let $u_j$ be
 smooth sections to $E$. This means that $u_j$ are functions that
depend smoothly on $t$ and holomorphically on $z$. To verify the
positivity of $E$ in the sense of Nakano we need to estimate from
below the
curvature of $E$ acting on the $k$-tuple $u$,
$$
\sum(\Theta^E_{j k}u_j,u_k).
$$
By (2.4) this means that we need to estimate from above
$$
\sum (\pi_\bot(\phi_ju_j),\pi_\bot(\phi_k u_k))=\|\pi_\bot(\sum \phi_ju_j)\|^2.
$$
Put $w=\pi_\bot(\sum \phi_ju_j)$.  For fixed $t$, $w$ solves the
$\dbar_z$-equation 
$$
\dbar w=\sum u_j\phi_{j \lambda}d\bar z_\lambda,
$$
since the $u_j$s  are holomorphic in $z$. Moreover, since $w$ lies in
the orthogonal complement of $A^2$, $w$ is the minimal
solution to this equation. 

We shall next apply Hörmander's weighted $L^2$-estimates for the
$\dbar$-equation. The precise form of these estimates that we need
says that if $f$ is a $\dbar$-closed form in a psedudoconvex domain 
$\Omega$, and if $\psi$ is a smooth strictly plurisubharmonic weight
function, then the minimal solution $w$ to the equation $\dbar v=f$
satisfies
$$
\int_\Omega |w|^2 e^{-\psi}\leq \int_\Omega\sum\psi^{\lambda
  \mu}f_\lambda\bar f_\mu e^{-\psi},
$$
where $(\psi^{\lambda \mu})$ is the inverse of the complex Hessian of
$\psi$ (see \cite{Demailly 1}).

In our case this means that 
$$
\int_\Omega|w|^2 e^{-\phi^t}\leq\int_\Omega\sum \phi^{\lambda
  \mu}\phi_{j \lambda}u_j\overline{\phi_{k\mu}u_k} e^{-\phi^t}.
$$
Inserting this estimate in formula (2.4) together with the formula for
the curvature of $F$ we find
\begin{equation}
\sum(\Theta^E_{j k}u_j,u_k)\geq\int_\Omega \sum_{j k}\left(\phi_{j
    k}-\sum_{\lambda \mu}\phi^{\lambda \mu}\phi_{j
    \lambda}\bar\phi_{k\mu}\right) u_j\bar u_k e^{-\phi^t}.
\end{equation}
We claim that the expression
$$
D_{j k}=:\left(\phi_{j
    k}-\sum_{\lambda \mu}\phi^{\lambda \mu}\phi_{j
    \lambda}\bar\phi_{k\mu}\right),
$$
in the integrand is a positive definite matrix at any fixed point . In
the proof of this, we may by 
a linear change
of variables in $t$  of course assume that the vector $u$ that
$D$ acts on equals $(1,0...0)$. The positivity of $D$ then follows
from a computation in \cite{Semmes}, but we will give the short
argument here too. Let $\Phi=i\ddbar \phi$ where the
$\ddbar$-operator acts on $t_1$ and the $z$-variables, the remaining
$t$-variables being fixed. Then
$$
\Phi=\Phi_{ 1 1} +i\alpha\wedge d\bar t_1+idt_1\wedge\bar\alpha +\Phi',
$$
where $\Phi_{1 1}$ is of bidegree $(1,1)$ in $t_1$, $\alpha$ is of
bidegree $(1,0)$  in $z$, and $\Phi'$ is of bidegree
$(1,1)$ in $z$. Then
$$
\Phi_{n+1}=\Phi^{n+1}/(n+1)!=
\Phi_{1,1}\wedge\Phi'_n-i\alpha\wedge\bar\alpha\wedge\Phi'_{n-1}\wedge
idt_1\wedge d\bar t_1.
$$
Both sides of this equation are forms of maximal degree that can be
written as certain coefficients multiplied by  the Euclidean volume form of
$\C^{n+1}$. The coefficient of the left hand side is the determinant
of the complex Hessian of
$\phi$ with respect to $t_1$ and $z$ together. Similarily, the
coefficient of the first term on the right hand side is $\phi_{1 1}$
times the Hessian of
$\phi$ with respect to the $z$-variables only. Finally, the
coefficient of the last term on the right hand side is the norm of the
$(0,1)$ form in $z$
$$
\dbar_z \partial_{t_1}\phi
$$
measured in the metric defined by $\Phi'$, multiplied by the volume
form of the same metric.
Dividing by the coefficient of $\Phi'_n$ we thus see that the matrix $D$
acting on a vector $u$ as above equals the Hessian of $\phi$ with
respect to $t_1$ and $z$ divided by the Hessian   of $\phi$ with respect to
the $z$-variables only. This expression is therefore positive so the
proof of Theorem 1.1 is complete. 

\section{Kähler fibrations with compact fibers} 

Let $X$ be a Kähler manifold of dimension $m+n$, fibered over a
complex $m$-dimensional manifold $Y$. This means that we have a
holomorphic map $p$ from $X$ to $Y$ with surjective differential at
all points. All our computations will be local, so we may as well
assume that $Y=U$ is a ball or polydisk in $\C^m$. For each $t$ in $U$
we let
$$
X_t=p^{-1}(t)
$$
be the fiber of $X$ over $t$. We shall assume that all fibers are
compact.

Next, we let $L$ be a holomorphic hermitian line bundle over $X$. Our
standing assumption on $L$ is that it is semipositive, i e that it is
equipped with a smooth hermitian metric of nonnegative curvature. For
each fiber $X_t$ we are interested in the space of holomorphic $L$-valued
$(n,0)$-forms on $X_t$,
$$
\Gamma(X_t, L|_{X_t}\gr K_{X_t})=: E_t.
$$
For each $t$, $E_t$ is a finite dimensional vector space and we claim
that
$$
E:=\bigcup\{t\}\times E_t
$$
has a natural structure as a holomorphic vector bundles.

To see this we need to study how $E_t$ varies with $t$. 
First note that $K_{X_t}$ is isomorphic to $K_X|_{X_t}$,
the
restriction of the canonical bundle of the total space to $X_t$, via
the map that sends a section $u$ to $K_{X_t}$ to
$$
\tilde u:= u\wedge dt,
$$
where $dt=dt_1\wedge ...dt_m$. It is clear that this map is
injective. Conversely, any local section $\tilde u$ to $K_X$ can be
locally represented as $\tilde u:= u\wedge dt$, and even though $u$ is
not uniquely determined, the restriction of $u$ to each fiber is
uniquely determined.
We thus have two ways of thinking of an element $u$ in $E_t$: as a holomorphic
$L$-valued $(n,0)$-form on $X_t$ or as a section $\tilde u=u\wedge dt$
to $K_X$ over
$X_t$. It turns out to be convenient for the computations later on to
have yet another 
interpretation: as the restriction of an $(n,0)$, $u'$, on $X$ to $X_t$ (here
we understand by restriction the pullback to $X_t$ under the inclusion
map from $X_t$ to $X$). Clearly, $u'$ is not uniquely determined by
$u$. Indeed $u'$ restricts to 0 on $X_{t_0}$ precisely when $u'\wedge dt=0$
vanishes for $t=t_0$, which in turn is equivalent to saying that
$$
u'=\sum\gamma_j\wedge dt_j.
$$
We will refer to a choice of $u'$ as a
representative of $u$. When $t$ varies, a smooth section of $E$ is
then represented by a smooth $(n,0)$-form on $X$. To avoid too
cumbersome notation we will in the 
sequel use the same letter to denote an element in $E_t$ and any
representative of it.

The semipositivity of $L$, and the assumption that $X$ is Kähler,
implies that any holomorphic section  $u$ to $K_{X_t}$ for one fixed $t$ can be
locally extended in the sense that there is a holomorphic section
$\tilde u$ to $K_X$ over $p^{-1}(W)$ for some neighbourhood $W$ of $t$
whose restriction to $X_t$ maps to $u$ under the isomorphism above. In
case $L$ is trivial this follows from the fact that Hodge numbers are
locally constant, see \cite{Voisin}. For general semipositive bundles
$L$ it follows from  a result of Ohsawa-Takegoshi type, that will be
discussed in an appendix. 

Taking a basis for $E_t$ for one fixed $t$ and extending as above we
therefore get a local frame for the bundle $E$. We define a complex
structure on $E$ by saying that an $(n,0)$-form over $p^{-1}(W)$,
$u$, whose restriction to each fiber is holomorphic, defines a holomorphic section of $E$
if $u\wedge dt$ is a holomorphic section of $K_X$. The frame we have
constructed is therefore holomorphic. 

Note that this means that $u$ is holomorphic if and only if $\dbar
u\wedge dt=0$, which means that  $\dbar
u$ can be written
\be
\dbar u=\sum\eta^j\wedge dt_j,
\ee
with $\eta^j$ smooth forms of bidegree $(n-1,1)$. Again, the $\eta^j$
are not uniquely determined, but their restrictions to fibers are. 

\smallskip

\noindent {\bf Remark:}
Even though we will not use it, it is worth mentioning the connection
between the forms $\eta_j$ and the Kodaira-Spencer map of the
fibration, see \cite{Voisin}, \cite{Kulikov}.

The Kodaira-Spencer map at a point $t$
in the base, is a map from the holomorphic tangent space of $U$ to the
first Dolbeault cohomology group, 
$$
H^{0,1}(X_t, T^{1,0}(X_t)),
$$
  of $X_t$ with values in the
holomorphic tangent space of $X_t$, i e, as $t$ varies  it is a
  $(1,0)$-form, $\sum \theta_j dt_j$  on $U$ with values  
in $H^{0,1}(X_t,T^{1,0}(X_t))$.The classes $\theta_j$  can be  represented 
by  $\dbar$-closed $(0,1)$-forms, $\vartheta_j$, on $X_t$ whose
coefficients are vector fields of type $(1,0)$ tangent to the
fiber. 
Such representatives can be found as follows.  Let $V_j$ be some
choice of smooth
$(1,0)$ vector fields on $X$, such that $dp(V_j)=\partial/\partial
t_j$. Then $dp(\dbar V_j)=0$ so $\dbar V_j$ are
$\dbar$-closed forms
with values in the bundle of vectors tangent to fibers. It is not hard
to check that they represent the classes $\theta_j$, e g using the
definition in \cite{Kulikov}.

Letting the 
vectorfield in the coefficients of $\vartheta_j$ act on forms by
contraction we  obtain maps 
$$
v\mapsto\vartheta_j\rfloor v
$$
from
$(p,q)$-forms on $X_t$ to $(p-1,q+1)$-forms. We claim that the forms
$$
\eta^j,
$$
restricted to fibers $X_t$,
is what we obtain when we let these map operate on $u$. Different
representatives of $u$ correspond to different
representatives of the same cohomology class.

To prove the claim, we need to verify that $\eta_j=\vartheta_j\rfloor
u$ on each fiber, where  $\vartheta_j$ is some representative on $X_t$
of the class $\theta_j$. 
Let  $ \widehat{dt_j}$ be the wedge product of all differentials
$dt_k$, except $dt_j$, with the right ordering, and let 
 $\tilde u=u\wedge dt$. Then
$$
V_j\rfloor\tilde u=(V_j\rfloor u)\wedge dt + u\wedge \widehat {dt_j}.
$$
Hence
$$
\dbar V_j\rfloor\tilde u=(\dbar(V_j\rfloor u))\wedge dt + \eta_j\wedge
dt.
$$
Since $\dbar V_j\rfloor dt=\dbar(V_j\rfloor dt)=0$, it follows that
$$
(\dbar V_j\rfloor u)\wedge dt= (\dbar(V_j\rfloor u))\wedge dt
+\eta_j\wedge dt.
$$
Therefore
$$
\dbar V_j\rfloor u=\dbar(V_j\rfloor u)+\eta_j
$$
on fibers, which proves our claim.
\qed

\smallskip

Let now $u$ be a smooth local section of $E$. This means that $u$ can
be represented by 
 a smooth $L$-valued form of bidegree $(n,0)$ over $p^{-1}(W)$ for some $W$
open in $U$, such that the restriction of $u$ to each fiber is
holomorphic. Then $\dbar u\wedge dt\wedge d\bar t=0$, so
$$
\dbar u =\sum d\bar t_j\wedge\nu^j +\sum \eta_j\wedge dt_j,
$$
where $\nu^j$ define sections to $E$. We define the $(0,1)$-part of
the connection $D$ on $E$ by letting
$$
D^{0,1}u=\sum \nu^j d\bar t_j.
$$
Sometimes we write 
$$
\nu^j=\dbar_{ t_j}u
$$
with the understanding that this refers to the $\dbar$ operator on
$E$. 
Note that $D^{0,1} u=0$ for $t=t_0$ if and only if each $\nu^j$ vanishes
when restricted to $X_{t_0}$, i e if $\dbar u\wedge dt=0$, which is
consistent with the definition of holomorphicity given earlier. Note
also that if we choose another $(n,0)$-form $u'$ to represent the same
section of $E$, then $u-u'$ vanishes when restricted to each
fiber. Hence $u-u'=\sum a_j\wedge dt_j$ and it follows that $D^{0,1}$
is well defined. 

The bundle $E$ has a naturally defined hermitian metric, induced by
the metric on $L$. To define the metric, let $u_t$ be an element of
$E_t$. Locally, with respect to a local trivialization of $L$, $u_t$
is given by a scalar valued $(n,0)$-form, $u'$, and the metric on $L$ is
given by a smooth weight function $\phi'$. Put
$$
[u_t,u_t]=c_n u'\wedge \overline{u'}e^{-\phi'},
$$
with $c_n=i^{n^2}$ is chosen to make this $(n,n)$-form
positive. Clearly this definition is independent of the
trivialization, so $[u_t,u_t]$ is globally defined. The metric on
$E_t$ is now defined as
$$
\|u_t\|^2=\int_{X_t}[u_t,u_t],
$$
and the associated scalar product is 
\be
(u_t,v_t)_t= \int_{X_t}[u_t,v_t].
\ee
In the sequel we will, abusively, write $[u,v]=c_nu\wedge\bar v e^{-\phi}$.
When $t$ varies we suppress the dependence on $t$ and get a smooth
hermitian metric on $E$. For local sections $u$ and $v$ to $E$ the
scalar product is then a function of $t$ and it will be convenient to
write this function as
$$
(u,v)=p_*([u,v])=p_*(c_nu\wedge\bar v e^{-\phi}),
$$
where $u$ and $v$ are forms on $X$ that represent the sections.
Here $p_*$ denotes the direct image, or push-forward, of a form,
defined by
$$
\int_U p_*(\alpha)\wedge\beta=\int_X \alpha\wedge p^*(\beta),
$$
if $\alpha$ is a form on $X$ and $\beta$ is a form on $U$. 

With the metric and the $\dbar$ operator defined on $E$ we can now
proceed to find the $(1,0)$-part of the Chern connection. Let $u$ be a
form on $X$ with values in $L$. Locally, with respect to a trivialization
of $L$, $u$ is given by a scalar valued form $u'$ and
the metric on $L$ is given by a function $\phi'$. Let
$$
\partial^{\phi'} u'=e^{\phi'}\partial (e^{-\phi'}u').
$$
One easily verifies that this expression is invariantly defined, and
we will, somewhat abusively, write $\partial^{\phi'} u'=\partial^{\phi}
u$, using $\phi$ to indicate the metric on $L$. Let now in particular
$u$ be of bidegree $(n,0)$ and such that the restrictions of $u$ to
fibers are holomorphic.
As $\partial^{\phi}u$ is of bidegree $(n+1,0)$ we can write
$$
\partial^{\phi}u=\sum dt_j\wedge \mu^j ,
$$
where $\mu^j$ are smooth $(n,0)$-forms whose restrictions to fibers are
uniquely defined. These restrictions are in general not holomorphic so
we let
$$
P(\mu^j)
$$
be the orthogonal projection of $\mu^j$ on the space of holomorphic
forms on each fiber. 
\begin{lma} The $(1,0)$-part of the Chern
connection on $E$ is given by
$$
D^{1,0}u =\sum P(\mu^j) dt_j. 
$$
\end{lma}

\begin {proof}
Even though it will follow implicitly from the proof below, we will
first prove that $D^{1,0}u$ is well defined, i e independent of the
choice of representative of $u$. Let therefore $u$ be a form that
restricts to 0 on $X_t$ for $t$ in some open set. Then we can write
$u=\sum\gamma_j\wedge dt_j$ there. Hence
$$
\partial^{\phi}u=\sum\partial^{\phi}\gamma_j\wedge dt_j,
$$
so $\mu_j=\partial^{\phi}\gamma_j$ on
 any fiber. It then follows from the definition of the
scalar product on $X_t$, and Stokes theorem, that $\mu_j$ is orthogonal
to all holomorphic forms on the fiber. In other words, $P(\mu_j)=0$,
so $D^{1,0}u=0$, which is what we wanted to prove.

To prove the lemma it suffices, by the definition of Chern connection, to
verify that
\be
\partial_{t_j}(u,v)=(P(\mu^j),v)+(u,\dbar_{t_j}v)=(\mu^j,v)+(u,\dbar_{t_j}v)
\ee
if $u$ and $v$ are smooth sections to $E$. But
$$
\partial(u,v)=\partial p_*([u,v])=
$$
$$
=c_n( p_*(\partial^\phi u\wedge
\bar v\,e^{-\phi})+(-1)^n p_*(u\wedge \overline{\dbar v}\,e^{-\phi}))=
$$
$$
c_n( p_*(\sum
\mu_j\wedge\bar v\wedge dt_j\, e^{-\phi}) +p_*(u\wedge\bar\nu^j\wedge
dt_j\, e^{-\phi})).
$$
This equals
$$
\sum ((\mu^j,v)+(u,\nu^j))dt_j,
$$
so we have proved 4.2.
\end{proof}

We will write $P(\mu^j)=D_{t_j}u$. 
We are now ready to verify the Nakano positivity of the bundle
$E$. For this we will use the recipe given at the end of section
2. Let $u_j$ be an $m$-tuple of holomorphic sections to  $E$ that
satisfy $D^{1,0} u_j=0$ at a given point that we take to be equal
to 0. Let 
$$
T_u=\sum (u_j,u_k)\widehat{dt_j\wedge d\bar t_k}.
$$
Here $\widehat{dt_j\wedge d\bar t_k}$ denotes the product of all
differentials $dt_i$ and $d\bar t_i$, except $dt_j$ and $d\bar t_k$
multiplied by a number of modulus 1, so that $T_u$ is nonnegative. We
need to verify that    
$$
i\ddbar T_u
$$
is negative. Represent the $u_j$s by smooth forms on $X$, and put
$$
\hat u=\sum u_j\wedge\widehat{dt_j}.
$$
Then, with $N=n+m-1$,
$$
T_u=c_Np_*(\hat u\wedge\overline{\hat u}\,e^{-\phi})
$$
Thus 
$$
\dbar T_u=c_N(p_{*}(\dbar\hat u\wedge\overline{\hat u} e^{-\phi})
+(-1)^{N}p_*(\hat u\wedge \overline{\partial^{\phi}\hat u}
e^{-\phi})). 
$$
Since each $u_j$ is holomorphic we have seen that
$$
\dbar u_j=\sum\eta^l_j\wedge dt_l.
$$
Therefore each term in the form
$$
p_{*}(\dbar\hat u\wedge\overline{\hat u} e^{-\phi})
$$
contains a factor $dt$. On the other hand, the push forward of an
$(n+m-1,n+m)$-form is of bidegree $(m-1,m)$, so we conclude that
$$
p_{*}(\dbar\hat u\wedge\overline{\hat u} e^{-\phi})=0.
$$
Thus
$$
\ddbar T_u=c_N((-1)^Np_*(\partial^\phi \hat
u\wedge\overline{\partial^\phi\hat u} e^{-\phi})+p_*(\hat
u\wedge\overline{\dbar\partial^\phi \hat u}))
$$
We rewrite the last term, using
$$
\dbar\partial^\phi+\partial^\phi\dbar=\ddbar\phi.
$$
Since
$$
p_*(\hat u\wedge\overline{\dbar\hat u}e^{-\phi})
$$
vanishes identically we find that
$$
(-1)^N p_*(\hat u\wedge\overline{\partial^\phi\dbar\hat
  u}e^{-\phi})+p_*(\dbar\hat u\wedge \overline{\dbar\hat u}e^{-\phi})=0,
$$
so all in all
\be
\ddbar T_u=
\ee
$$
c_N\left((-1)^N p_*(\partial^\phi \hat
u\wedge\overline{\partial^\phi\hat u} e^{-\phi})-p_*(\hat
u\wedge\overline{\hat u}\wedge\ddbar\phi e^{-\phi})+(-1)^N p_*(\dbar\hat
u\wedge\overline{\dbar\hat u}e^{-\phi})\right).
$$
So far, the computations hold for any choice of representative of our
sections $u_j$. We shall next choose our representatives in a careful way.
\begin{prop}
Let $u$ be a section of $E$ over an open set $U$ containing
the origin, such that
$$
D^{0,1} u=0,
$$
in $U$ and
$$
D^{1,0}u=0
$$
at $t=0$. Then u can be represented by a smooth $(n,0)$-form, still
denoted $u$, such that
\be
\dbar u =\sum\eta^k\wedge dt_k,
\ee
where $\eta^k$ is {\bf primitive} on $X_0$, i e satisfies
$\eta^k\wedge \omega=0$ on $X_0$, and furthermore
\be
\partial^{\phi}u\wedge \widehat dt_j=0,
\ee
at $t=0$ for all $j$.
\end{prop}

To prove the proposition we need two lemmas.
\begin{lma} Let $u$ be an $(n,0)$-form on $X$, representing a
  holomorphic section of $E$, and write
$$
\dbar u=\sum\eta^k\wedge dt_k.
$$
Then $\eta^k\wedge \omega$ are $\dbar$-exact on any fiber.
\end{lma}

\begin{proof}
Since $u\wedge\omega$ is of bidegree $(n+1,1)$ we can write locally
$$
u\wedge\omega=\sum u^k\wedge dt_k.
$$
The coefficients $u^k$ here are not unique, but their restrictions to
fibers are unique. This follows since $\sum u^k\wedge dt_k=0$ implies
$$
\sum u^k\wedge dt_k\wedge \widehat dt_k=u^k\wedge dt=0,
$$
which implies that $u^k$ vanishes when restricted to any fiber.

Hence in particular $u^k$ are well defined global forms on any fiber.
Moreover
$$
\sum\eta^k\wedge\omega\wedge dt_k=\dbar u\wedge\omega=\sum \dbar u^k\wedge dt_k,
$$
so
$$
\sum(\eta^k\wedge\omega-\dbar u^k)\wedge dt_k=0.
$$
Again, wedging with $\widehat dt_k$, we see that 
$$
\eta^k\wedge\omega=\dbar u^k
$$
on fibers, so $\eta^k\wedge\omega$ is exact on fibers.
\end{proof}

\begin{lma}
Let $\mu$ be an $(n,0)$-form on a compact $n$-dimensional Kähler manifold
$Z$, with values in a hermitian holomorphic line bundle $L$. Assume $\mu$ is
orthogonal to the space of holomorphic $L$-valued forms under the
scalar product (4.2). Let $\xi$ be a $\dbar$-exact $(n,2)$-form on $Z$,
with values in $L$. Then there is an $L$-valued form $\gamma$
of bidegree $(n-1,0)$ such that
$$
\partial^{\phi}\gamma=\mu,
$$
and
$$
\dbar\gamma\wedge \omega=\xi.
$$
\end{lma}
\begin{proof}
Since $\xi$ is exact we can solve $\dbar\chi=\xi$. Then
$\mu-\dbar^*\chi$ is orthogonal to holomorphic forms, so we can solve
$$
\dbar^*\alpha=\mu-\dbar^*\chi,
$$
with $\dbar\alpha=0$. This follows since the range of $\dbar^*$ is
closed on a compact manifold. Let $\alpha'=\alpha+\chi$, so
$$
\mu=\dbar^*\alpha'
$$
and $\dbar\alpha'=\xi$. Write $\alpha'=\gamma\wedge\omega$, where
$\omega$ is the Kahler form. Then $\gamma$ satisfies the conditions in
the lemma (possibly up to a sign).
\end{proof}

We are now ready to prove the proposition  4.2.

Recall that $u$ can in any case be represented by a form satisfying
$$
\dbar u=\sum \eta^k\wedge dt_k,
$$
over $U$ and
$$
\partial^\phi u=\sum \mu^k\wedge dt_k,
$$
where the restriction of $\mu^k$ is orthogonal to holomorphic forms on
$X_0$. By Lemmas 4.3 and 4.4  there are forms $\gamma^k$ on $X_0$ such
that
$$
\mu^k=\partial^\phi\gamma^k
$$
and
$$
\dbar\gamma^k\wedge\omega=\eta^k\wedge\omega
$$
on $X_0$. Extend $\gamma^k$ smoothly to a neighbourhood of $X_0$ ( i e
find a form that restricts to $\gamma^k$), and put
$$
u'=u-\sum\gamma^k\wedge dt_k.
$$
Then $u'$ is a form of bidegree $(n,0)$ that  represents the same section of
$E$ as $u$.
Then
$$
\partial^\phi u'\wedge \widehat dt_j=\partial^\phi u\wedge \widehat dt_j
-\partial^\phi\gamma^j\wedge dt=(\mu_j-\partial^\phi\gamma^j)\wedge
dt=0
$$
{\it at} $t=0$, since 
$$
\mu_j=\partial^\phi\gamma^j
$$
{\it on} $X_0$. Moreover
$$
\dbar u'=\sum(\eta^k-\dbar\gamma^k)\wedge dt_k,
$$
and
$$
((\eta^k-\dbar\gamma^k)\wedge \omega=0
$$
on $X_0$. Hence $u'$ satisfies all the requirements and the
propositions is proved.

\bigskip

We now return to the proof of Theorem 1.2. Note that with the choice
of representatives of our sections $u_j$ furnished by Proposition 4.2,
formula 4.4 
simplifies at $t=0$ to
\be
\ddbar T_u=
c_N\left(-p_*(\hat
u\wedge\overline{\hat u}\wedge\ddbar\phi e^{-\phi})+(-1)^N p_*(\dbar\hat
u\wedge\overline{\dbar\hat u}e^{-\phi})\right).
\ee
The first term on the right hand side obviously gives a (semi)negative
contribution to $i\ddbar T$. 
To analyse the last term write
$$
\dbar u_j=\sum \eta^k_j\wedge dt_k,
$$
and
$$
\dbar\hat u=\sum \eta^j_j \wedge dt=:\eta\wedge dt.
$$
Then the last term equals
$$
c_n\int_{X_0}\eta\wedge\bar\eta \, \,e^{-\phi}\,dV_t.
$$
In general the quadratic form in $\eta$ appearing here is
indefinite. In our case however, all the $\eta^k_j$ are primitive on
$X_0$, and it is well known that
$$
c_n \eta\wedge\bar\eta=-|\eta|^2 
$$
if $\eta$ is primitive. (This is easily checked by hand at a point by
choosing coordinates that are orthogonal at the point.) Hence
$$
c_n\int_{X_0}\eta\wedge\bar\eta e^{-\phi}\,dV_t=-\int_{X_0}|\eta|^2 \,dV_t.
$$
and we get

\be
i\ddbar T= c_N\left(-p_*(\hat u\wedge\overline{\hat u}\wedge
  i\ddbar\phi e^{-\phi})\right) 
-\int_{X_0}|\eta|^2 \,dV_t.
\ee

By 4.4, this means that $i\ddbar T_u\leq 0$, so $E$ is at least
seminegative in the sense of Nakano. If $i\ddbar\phi$ is strictly
positive, it is clear that the curvature term alone gives a strictly
negative contribution to $i\ddbar T$.  Therefore $E$ is
strictly positive if 
$L$ is strictly positive so we have proved Theorem 1.2.
In   the next section we shall see that even when $L$ is only
semipositive, equality can hold in our estimates only in very special
cases.

We want to add one  remark on the relation between the proof of
Theorem 1.2 in this section and the proof of Theorem 1.1 in section 3.
The proof in section 3 is easily adapted to the case of a trivial
fibration (so that $X$ is a global product). It may then seem that the
proof here is quite different since it does not use the
Hörmander-Kodaira $L^2$-estimates at all.
The two proofs are however really quite similar, the difference being
that in this section we basically reprove the special case of the
$L^2$-estimates that we need as we go along.

\section{Semipositive vector bundles.}
In this section we will discuss when equality holds in the inequalities
of Theorem 1.2, i e when the bundle $E$ is not strictly positive. As
we have already seen in the last section, this can only happen if the
line bundle $L$ is not strictly positive. More precisely, provided the
components $u_j$ of $\hat u$ are chosen to satisfy the conditions in
Proposition 4.2,  equality
holds if and only if
$$
\eta=0
$$
and
$$
\hat u\wedge\bar{\hat u}\wedge i\ddbar\phi =0,
$$
and since $i\ddbar \phi\geq 0$, the last condition is equivalent to
\be
\hat u\wedge i\ddbar\phi =0.
\ee

For simplicity we assume from now that the base domain $U$ is
one-dimensional, so that we do not need to discuss degeneracy in
different directions, and 
that the fibration we consider is locally holomorphically trivial, i e
that $X=U\times Z$, where $Z$ is a compact $n$-dimensional complex
manifold. 
Moreover, we assume that the curvature of our metric $\phi$ on $L$,
$\Theta^L$, 
is strictly positive along each fiber $X_t\simeq Z$. Then we can also
assume that we have chosen our Kähler metric $\omega$ on $X$ so that
$\omega_t:=\omega|X_t= i\Theta^L|X_t$ on each fiber. 

Since $X$ now is a global product we can decompose $\Theta^L$
according to its degree in $t$ and $z$, where $z$ is a 
ny local coordinate on $Z$. In particular, there is a well defined
$(0,1)$-form $\theta^L$ on $X$ such that $dt\wedge\theta^L$ is the
component of $\Theta^L$ of degree 1 in $dt$. Expressed in invariant
language,
$$
\theta^L=\delta_{\partial/\partial t} \Theta^L,
$$
where $\delta_{\cdot}$ means contraction with a vector field. Using
the formulas
\be
\partial\delta_V +\delta_V \partial=\mathcal L_V,
\ee
 where $\mathcal L$ is the holomorphic Lie derivative, and
\be
\dbar\delta_V+\delta_V\dbar =0,
\ee
if $V$ is a holomorhic vector field we see that $\theta^L$ is
$\dbar$-closed on $X$, and that
$$
\partial \theta^L=\mathcal L_{\partial/\partial t}
\Theta^L
$$
on $X$.
On the other hand, on each fiber $X_t$ there is a 
unique $(1,0)$ vector field $V_t$, defined by
$$
\delta_{V_t}\omega_t=\theta^L.
$$
Our key obeservation is contained in the next lemma. 
\begin{lma}
Assume that for some $u\neq 0$ in $E_0$, $(\Theta^E u,u)=0$. Then
$V_0$ is a holomorphic vector field on $X_0$.
\end{lma}
\begin{proof}

Since the base $U$ is one dimensional, $\hat u=u$ and since $X=Z\times
U$ is a global product we can decompose
$$
u=u_0 + dt\wedge v,
$$
where $u_0$ does not contain $dt$, i e $\delta_{\partial/\partial
  t}u_0=0$. Here $u$ is chosen to satisfy the conditions of
proposition 4.2,
$$
\dbar u=\eta\wedge dt,
$$
where $\eta\wedge \omega=0$ on $X_0$, and
$$
\partial^\phi u=0
$$
for $t=0$. Since $u$ is holomorphic on $X_0$, $\dbar u_0=0$, so
$\eta=\dbar v$ on $X_0$.
But we have seen above that if $u$ is a null vector for the
curvature of $E$ at $t=0$, then $\eta=0$ so
the restriction of $v$ to $X_0$
is a holomorphic $(n-1,0)$-form. 

We also know from (5.1) that
$$
u\wedge\Theta^L=0
$$
for $t=0$. Applying $\delta_{\partial/\partial t}$, we get that
$$
v\wedge\Theta^L + (-1)^nu\wedge \theta^L=0
$$
for $t=0$. Restricting to $X_0$
$$
v\wedge \Theta^L=(-1)^{n+1}u\wedge\delta_{V_0}\Theta^L
$$
on $X_0$. But, $u\wedge\Theta^L=0$ on $X_0$ (for reasons of bidegree)
so
$$
\delta_{V_0}u\wedge\Theta^L +(-1)^n u\wedge\delta_{V_0}\Theta^L=0.
$$
Hence $v$ equals $\delta_{V_0}u$ on $X_0$, so
$\delta_{V_0}u$ is a holomorphic form on $X_0$. Since $u$ is also
holomorphic it follows that $V_0$ 
must be holomorphic too, except possibly where $u$ vanishes. But since
$V$ is smooth, $V$ must actually be holomorphic everywhere by
Riemann's theorem on removable singularities.
\end{proof}

This lemma could also have been proved using the approach via
$L^2$-estimates. It is then strongly related to the following
proposition that we state explicitly since we feel it has an
independent interest.

\begin{prop}Let $L$ be a positive line bundle over a compact complex
  manifold $Z$. Give $Z$ the Kähler metric defined by the curvature
  form of $L$. Let $\mu$ be the $L^2$-minimal solution to
  $\dbar\mu=f$, where $f$ is an $L$-valued $(n,1)$-form on $Z$. Then
  equality holds in Hörmander's estimate, i e
$$
\int_Z|\mu|^2=\int_Z |f|^2
$$
if and only if $\gamma=*f$ is a holomorphic form.
\end{prop}
\begin{proof}
Let $\phi$ be the metric on $L$. By Lemma 4.4
$$
\mu=\partial^\phi\gamma,
$$
for some $\dbar$-closed $(n-1,0)$-form $\gamma$. Thus
$$
\int_{Z}|\mu|^2=\int_Z f\wedge\bar\gamma e^{-\phi}\leq\|f\|\|\gamma\|,
$$
with equality only if $*f$ is proportional to $\gamma$. 
By the
Hörmander-Kodaira-Nakano identity
$$
\int_{Z}\gamma\wedge\bar\gamma\wedge
i\ddbar\phi\,e^{-\phi}+\int_{Z}|\dbar\gamma|^2=\int_{Z}|\mu|^2.
$$
The first term on left hand side here is the norm squared of $\gamma$
so it follows 
that
$$
\|\gamma\|^2\leq \|\mu\|^2,
$$
with equality only if $\dbar \gamma=0$, and combined with our previous estimate
$$
\|\mu\|^2\leq\|f\|^2
$$
with equality only if $\dbar\gamma=0$ and $*f$ is proportional to
$\gamma$. Hence $*f$ must be holomorphic.  The argument is easily seen
to be reversible. 
\end{proof}

We are now ready to state the main theorem of this section.
\begin{thm}Assume that $Z$ has no nonzero global holomorphic vector
  field. Suppose  that

(i)$X$ is locally a product $U\times Z$ where $U$
  is an open set in $\C$,

(ii) $L$ is semipositive on $X$,

and that 

(iii) $L$ restricted to each fiber
  is strictly positive. 
Let $\omega_t$ the Kähler metric be the
fiber $X_t$ induced by the curvature of $L$. Then, if for
each $t$ in $U$ there is some element $u_t$ in $E_t$ such that
$$
(\Theta^E u_t,u_t)=0,
$$
it follows that
$$
\omega_t=\omega_0
$$
for $t$ in $U$.
\end{thm}
\begin{proof}

\noindent By Lemma 5.1 the restriction of $\theta^L$ to each fiber
$X_t$ is zero. Hence 
$$
\partial\theta^L=\mathcal L_{\partial/\partial t}\Theta^L
$$
also vanishes on fibers, which means that
$$
\frac{d}{dt}\omega_t=0.
$$
\end{proof}

\section{The space of Kähler metrics.}

In this section we will specify the situation even more, and assume
that $X=U\times Z$ is a product, and that moreover the line bundle $L$
is the pullback of a bundle on $Z$ under the projection on the second
factor. Intuitively this means that not only  are all fibers 
the same, but also the line bundle on them, so it is only the metric
that varies. Fix one metric $\phi_0$ on
$L$, that we can take to be the pullback of a metric on the bundle on
$Z$, i e independent of the $t$-variable. Then any other metric on $L$
can be written
$$
\phi=\phi_0+\psi,
$$
where $\psi$ is a function on $X$. 
We also continue to assume
that $U$ is a domain in $\C$. Let $u$ be an element in $E_t$.

In this situation we have an explicit lower bound for the curvature
form  operating on $u$, generalizing 3.1:
\be
(\Theta^E u,u)\geq \int_{X_t}\left(\psi_{t \bar
  t}-|\dbar_z\psi_t|^2_{\phi}\right)[u,u].
\ee
Here the expression $|f|_\phi$ means the norm of the form $f$ with
respect to the metric $\omega:=i\ddbar_z\phi$ on $X_t$. 
This can be proved, either by adapting the method of section 3 - note
that we may replace any $t$-derivative of $\phi$ by the corresponding
derivative of $\psi$ since $\phi_0$ is independent of $t$ - or from 
 the more complicated proof in section 4. To see how it follows from
 the formulas in section 4 we again decompose
$$
u=u_0 +dt\wedge v
$$
like in the previous section, and also write
$$
i\ddbar\phi= \omega -2\Re idt\wedge\dbar\psi_t +\psi_{t\bar
  t}idt\wedge d\bar t.
$$
Then
$$
c_N u\wedge\bar u\wedge i\ddbar\phi=\left(\psi_{t \bar t}c_n u_0\wedge\bar u_0
+c_{n-1}v\wedge\bar v\wedge\omega-2\Re\dbar\psi_t\wedge u_0\wedge \bar
v\right)\wedge idt\wedge d\bar t.
$$
By Cauchy's inequality
$$
2\Re\dbar\psi_t\wedge u_0\wedge \bar v\leq c_{n-1}v\wedge\bar v\wedge
\omega +|\dbar\psi_t|_\phi^2c_n u_0\wedge\bar u_0
$$
so
$$
c_N u\wedge\bar u\wedge i\ddbar\phi\geq\left( \psi_{t\bar
    t}-|\dbar\psi_t|_\phi^2\right)c_n u_0\wedge\bar u_0,
$$
and (6.1) follows from (4.8)  and (2.5). Notice that we have used
nowhere that $\Theta$ is positive on the total space $X$, just that
the restriction to fibers are positive. Therefore (6.1) holds for any
metric on $L$ which is strictly positive along the fibers, even though
of course this does not imply that $E$ is positive in general.

The expression occuring in the integrand in 6.1,
$$
C(\psi)= \left(\psi_{t \bar t}-|\dbar_z\psi_t|^2_{\phi}\right)
$$
plays a crucial role
in the recent work on variations of Kähler metrics on compact
manifolds, see \cite{Semmes}, \cite{Mabuchi}, \cite{Donaldson1},
\cite{Donaldson2},\cite{Phong-Sturm} and \cite{Chen-Tian}, to quote
just a few. Fixing a line bundle $L$ on $Z$, these papers consider the
space $\mathcal K(L)$ of all Kähler metrics whose Kähler form is
cohomologous to the Chern class of $L$. This means precisely that the
Kähler form can be written
$$
i\ddbar\phi=i\ddbar\phi_0+i\ddbar\psi,
$$
for some function $\psi$. and so the set up we described above, where
$\psi$ depends on $t$, corresponds to a path in $\mathcal K(L)$.

The tangent space of $\mathcal K(L)$ at a point $\phi$ is a space of functions
$\dot{\psi}$ and a Riemannian metric on the tangent space is given by
the $L^2$-norm
$$
|\dot{\psi}|^2=\int_Z|\dot{\psi}|^2 (i\ddbar\phi)^n/n!.
$$
In this way, $\mathcal K(L)$ becomes an infinite-dimensional
Riemannian manifold.

Now consider our space $X$ above and let $U=\{|\Re t|<1\}$ be a strip,
and consider functions $\psi$ that depend only on $\Re t$. Then
$$
4C(\psi)=\ddot{\psi}- |\dbar_z\dot{\psi}|^2_\phi
$$
if we use dots to denote derivatives with respect to $\Re t$. The link
between Theorem 1.2 and the papers quoted above lies in the fact that,
by the results in \cite{Semmes}, the right hand side here is the
geodesic curvature of the path in $\mathcal K(L)$ determined by
$\psi$. 

A basic idea in the papers quoted above is to consider the spaces
$$
E_t=\Gamma(X_t, K_{X_t}\gr L)
$$
with the induced $L^2$-metric
as a finite dimensional approximation or ``quantization'' of the
manifold $Z$ with metric $\phi_t=\phi_0 +\psi(t,\cdot)$. (Actually
this is not quite true. In the papers quoted above one does not take
the tensor product with the canonical bundle, but instead integrates
with respect to the volume element $(i\ddbar\phi_t)^n/n!$.) Here one also
replaces $L$ by $L^k$ - with $k^{-1}$ playing the role of Planck's
constant - and studies the asymptotic behaviour as $k$ goes to 
infinity. 

Under this ``quantization'' map, functions, $\chi$, on $Z$ correspond 
to the induced Toepliz operator, $T_{\chi}$,  on $E_t$. This Toepliz
operator is defined by
$$
(T_{\chi}u,u)_{E_t}=\int_{\{t\}\times Z}\chi [u,u],
$$
if $u$ is any element in $E_t$. 
Note that the right hand side in our estimate for the curvature (6.1)
equals ( 4 times) $(T_{\chi}u,u)$, with $\chi$ equal to the geodesic
curvature of the path in $\mathcal K$. 
Thus the inequality 6.1 can  be formulated as saying that ``the curvature
of the quantization is greater than the quantization of the
curvature'', i e that the curvature operator of the vector bundle corresponding
to a path in $\mathcal K(L)$ is greater than the Toepliz operator
defined by the geodesic curvature of the path. Moreover, Theorem 5.3
implies that if $Z$ has no nonzero global holomorphic vector fields,
then equality holds only for a constant path.

\section{ Bundles of projective spaces}
Let $V$ be a  holomorphic vector bundle of finite rank $r$ over a complex
manifold $Y$, and let $V^*$ be its dual bundle. We let $\mathbb P(V)$
be the fiber bundle over $Y$ whose fiber at each point $t$ of the base
is the projective space of lines in $V^*_t$, $\mathbb P(V^*_t)$. Then
$\mathbb P(V)$ is a holomorphically locally trivial fibration. There
is a naturally defined line bundle $O_{\mathbb P(V)}(1)$ over $\mathbb
P(V)$ whose restriction to any  fiber $\mathbb P(V^*_t)$ is the hyperplane
section bundle (see \cite{Lazarsfeld}). One way to define this bundle is to first consider the
tautological line bundle
$O_{\mathbb P(V)}(-1)$. The total space of this line bundle, with the
zero section removed, is just
the total space of $V^*$ with the zero section removed, and the
projection to  $\mathbb P(V)$ is the map that sends a nonzero point in
$V^*_t$ to its image in  $\mathbb P(V^*_t)$. The bundle 
$O_{\mathbb P(V)}(1)$ is then defined as the dual of $O_{\mathbb P(V)}(-1)$. 
The global holomorphic sections
of this bundle over any fiber are in one to one correspondence with
the linear forms on $V_t^*$, i e the elements of $V$. More generally, $
O_{\mathbb P(V)}(1)^l=O_{\mathbb P(V)}(l)$ has as global holomorphic
sections over each fiber the homogenuous polynomials on $V^*_t$ of
degree $l$, i e the elemets of the $l$:th symmetric power of $V$.
We shall apply Theorem 1.2 to the line bundles 
$$
L(l)=: O_{\mathbb P(V)}(l).
$$

Let  $E(l)$ be the vector bundle whose fiber over a point $t$ in $Y$ is
the space of global holomorphic sections of 
$L(l)\gr K_{\mathbb P(V^*_t)}$. If $l<r$ there is only the zero
section, so we assume 
from now on that $l$ is greater than or equal to $r$. 

We claim that
$$
E(r)=\det V,
$$
the determinant bundle of $V$. To see this, note that $L(r)\gr
K_{\mathbb P(V^*_t)}$ is trivial on each fiber, since the canonical
bundle of $(r-1)$-dimensional projective space is $O(-r)$. The space
of global sections is therefore one dimensional. A convenient  basis
element is
$$
\sum_1^r z_j \widehat{dz_j},
$$
if $z_j$ are coordinates on $V^*_t$. Here $\widehat{dz_j}$ is the
wedge product of all differentials $dz_k$ except $dz_j$ with a sign
chosen so that $dz_j\wedge\widehat{dz_j}= dz_1\wedge ...dz_r$. If we
make a linear change of 
coordinates on $V^*_t$, this basis element gets multiplied with the
determinant of the matrix giving the change of coordinates, so the
bundle of sections must transform as the determinant of $V$. Since
$$
L(r+1)\gr K_{\mathbb P(V^*_t)}=O_{\mathbb P(V)}(1)\gr L(r)\gr
K_{\mathbb P(V^*_t)},
$$
it also follows that
$$
E(r+1)=V\gr\det V.
$$ 
In the same way 
$$
E(r+m)= S^m(V)\gr\det V,
$$
where  $S^m(V)$ is the $m$th symmetric power of $V$. 

Let us now assume that $V$ is ample in the sense of Hartshorne, see
\cite{Harts}. By a theorem of Hartshorne, \cite{Harts}, $V$ is ample
if and only if 
 $L(1)$ is  ample, i e has a metric with
strictly positive curvature. Theorem 1.2 then implies that the
$L^2$-metric on each of the bundles $E(r+m)$
for $m\geq 0$ has curvature which is strictly positive in the sense of
Nakano, so we obtain:
\begin{thm} Let $V$ be a vector bundle (of finite rank) over a complex
  manifold. 
Assume $V$ is ample in the sense of Hartshorne. Then for any $m\geq 0$
the bundle
$$
S^m(V)\gr\det V
$$
has an hermitian metric with curvature which is (strictly) positive in
the sense of Nakano.
\end{thm}

\section{Appendix}
In this section we will state and prove an extension  result of Ohsawa-Takegoshi
type which in particular implies that the bundles $E$ that we have
discussed in this paper really are vector bundles. The proof follows
the method of \cite{Berndtsson 2}. See also \cite{Siu} for a closely
related result. 
\begin{thm}
Let $X$ be a Kähler manifold fibered over the unit ball $U$ in $\C^m$, with
compact fibers $X_t$. Let $L$ be a holomorphic line bundle on $X$ with
a smooth hermitian metric  with semipositive
curvature. Let $u$ be a holomorphic section of $K_{X_0}\gr L$ over
$X_0$ such that
$$
\int_{X_0}[u,u]\leq 1.
$$
Then there is a holomorphic section, $\tilde u$ to $K_X$ over $X$ such
that $\tilde u=u\wedge dt$ for $t=0$ and
$$
\int_X [\tilde u,\tilde u]\leq C
$$
where $C$ is an absolute constant.
\end{thm}
\begin{proof}
We assume $m=1$. The general case follows in the same way, extending
with respect to one variable at the time. At first we also assume that
the metric on  $L$ is smooth. The proof follows closely the method in
\cite{Berndtsson 2} so we will be somewhat sketchy.

Let $f=u\wedge [X_0]/(2\pi i)$, where $[X_0]$ is the current of integration on $X_0$.
Then $\dbar f=0$ and if $v$ is any solution to $\dbar v =f$ then
$\tilde u= tv$ is a section of $K_X\gr L$ that extends $u$ in the sense
described. To find a $v$  with $L^2$-estimates we need to estimate
$$
\int_X (f,\alpha)
$$
for any compactly supported test form $\alpha$ of bidegree $(n+1,1)$
on $X$. For $\alpha$ given, decompose $\alpha=\alpha^1 +\alpha^2$,
where $\alpha^1$ is $\dbar$-closed, and $\alpha^2$ is orthogonal to
the kernel of $\dbar$. This means that $\alpha^2$ can be written
$$
\alpha^2 =\dbar^* \beta
$$

for some $\beta$. By the regularity of the $\dbar$-Neumann problem
$\alpha^i$ and $\beta$ are all smooth up to the boundary. We first
claim that
\be
\int (f,\alpha^2) =0.
\ee
This is not surprising since $f$ is $\dbar$-closed, but it is not
quite evident since $f$ is not in $L^2$. To prove it, extend $u$
smoothly to $X$. Then $\dbar u\wedge d\bar t=0$ for $t=0$. Let $\chi$
be a smooth cut-off function equal to one near the origin in $\R$, and
put 
$$
\chi_\epsilon(t)=\chi(|t|^2/\epsilon).
$$
Then 
$$
f=\chi_\epsilon f=\dbar(u\wedge \frac{dt}{t})\chi_\epsilon-\dbar u\wedge
\frac{dt}{t}\chi_\epsilon=
$$
$$
=\dbar(u\wedge
\frac{dt}{t}\chi_\epsilon)- u\wedge \frac{dt}{t}\wedge\dbar \chi_\epsilon-\dbar u\wedge
\frac{dt}{t}\chi_\epsilon=: I + II + III.
$$
Clearly the scalar product between I and $\alpha^2$ vanishes. It is
also clear that the scalar product between III and $\alpha^2$ goes to
zero as $\epsilon$ goes to zero. The scalar product between II and
$\alpha^2$ equals, up to signs
$$
\int\chi' (\dbar u\wedge dt\wedge d\bar t,\beta )/\epsilon,
$$
which is easily seen to tend to zero as well since $\dbar u\wedge
d\bar t$ vanishes for $t=0$. Hence 8.1
follows. Therefore
$$
|\int_X (f,\alpha)|^2=|\int_X
(f,\alpha^1)|^2\leq\int_{X_0}\gamma\wedge\bar\gamma e^{-\phi},
$$
where $\gamma$ is the Hodge-* of $\alpha^1$.  The form $\gamma$
satifies $\omega\wedge\gamma=\alpha^1$ and
$\dbar^*\alpha=\dbar^*\alpha^1=\partial^{\phi}\gamma$. To estimate this we apply
the Siu $\ddbar$-Bochner formula ( see \cite{Berndtsson 2}): If $w$ is
any nonnegative function smooth up to the boundary of $X$, then
\be
-\int i\ddbar w\wedge c_n\gamma\wedge\bar\gamma e^{-\phi}+
\int i\ddbar \phi\wedge c_n\gamma\wedge\bar\gamma e^{-\phi}w\leq
\ee
$$
\leq
2c_n\Re\int \dbar\partial^\phi\gamma\wedge\bar\gamma e^{-\phi}w=
2\int |\dbar^*\alpha|^2w +2\Re\int (\dbar^*\alpha,\partial w\wedge\gamma).
$$
Now choose $w=(1/2\pi)\log(1/|t|^2)$. (Although $w$ is not smooth it
can be approximated by the smooth functions
$(1/2\pi)\log(1/(|t|^2+\epsilon))$, so formula (8.2) still holds.)  If
$i\ddbar\phi$ is nonnegative 
we then find that
$$
\int_{X_0}\gamma\wedge\bar\gamma e^{-\phi}\leq C\int |\dbar^*\alpha|^2
(\log(1/|t|^2+1/|t|) +C\int_{X}idt\wedge d\bar
t\wedge\gamma\wedge\bar\gamma e^{-\phi} (1/|t|).
$$
To take care of the last term we repeat the last argument once more,
this time choosing $w=(1-|t|)$ and finally obtain an estimate
$$
|\int_X (f,\alpha)|^2\leq C\int |\dbar^*\alpha|^2(1/|t|).
$$
This implies that there is some function $v$ on $X$ such that
$$
\int_X (f,\alpha)=\int(v,\dbar^*\alpha),
$$
for all test forms $\alpha$, and satisfying
$$
\int |v|^2 |t|\leq C.
$$
Then $\tilde u:=tv$ satisfies the conclusion of the theorem.

\end{proof}

\bigskip

\def\listing#1#2#3{{\sc #1}:\ {\it #2}, \ #3.}

\end{document}